\documentclass[11pt,twoside]{article}
\usepackage{amsmath,amssymb}
\usepackage{graphicx}
\usepackage[footnotesize]{subfigure}
\usepackage[footnotesize]{caption2}

\allowdisplaybreaks[4]

\newtheorem{lemma}{{\bf Lemma}}[section]

\newtheorem{theorem}{{\bf Theorem}}[section]

\makeatletter
\evensidemargin
\oddsidemargin
\makeatother

\title{
\Large\bf Global dynamics of a competition-diffusion system and
\\ application to a modified Leslie-Gower model}

\author{{\sc Leqi Chen$^{a}$}, {\sc Shuang Chen$^{b}$\footnote{\rm{Corresponding Author}: {schen@hust.edu.cn}}}
\\
{\small $^{a}$ School of Mathematics, Sichuan University}\\
{\small Chengdu, Sichuan 610064, P. R. China}\\
{\small $^{b}$ School of Mathematics and Statistics, Huazhong University of Sciences and Technology}\\
{\small Wuhan, Hubei 430074, P. R. China}
}

\date{}

\begin{document}
\maketitle
\begin{abstract}
We investigate the global dynamics of a  Lotka-Volterra competition-diffusion system
in spatially heterogeneous environment.
This model indicates that the evolution of the density of the predator is independent of the density of the prey.
Based on the  principal spectral theory and the dynamics of the classical single-species logistic model,
we obtain the global dynamics of this competition-diffusion system.
As an application, under some suitable conditions
we use the obtained results to prove the global stability of steady states and the persistence of the two species
in a modified Leslie-Gower model with diffusion in heterogeneous environment.

\vskip 0.2cm
{\bf AMS(2010) Subject Classification}: 35Q92, 92D25, 35K57, 35P05.
\end{abstract}
\baselineskip 15pt
\parskip 10pt

\thispagestyle{empty}
\setcounter{page}{1}


\section{Introduction}
\setcounter{equation}{0}
\setcounter{lemma}{0}
\setcounter{theorem}{0}
\setcounter{remark}{0}

The Lotka-Volterra models, which are originated by Lotka \cite{Lotka-1910} and Volterra \cite{Volterra-26,Volterra-31},
have been frequently used to describe population dynamics in the past couple of decades.
Among the investigation to the interactions between movement and environmental heterogeneity in population dynamics,
many efforts have been devoted to studying the dynamics of single species models \cite{DeAngelis-16,Liang-Lou-12,Lou-06},
Lotka-Volterra competition systems
\cite{Avila-Vales-97,Cantrelletal-91,He-Ni-13-1,He-Ni-13-2,He-Ni-16,Lam-12,Lam-15,Lou-19,Tang-Zhou-20,Zhou-Xiao-18},
and Lotka-Volterra predator-prey systems \cite{Du-Hsu-04,Du-Shi-07,Lou-17,peng-09,Wang-Li-13,Zhao-Wang-14,Zou-Guo}
with diffusion in spatially heterogeneous environment.
For more detail concerning  this topic,
we also refer the readers to excellent monographs \cite{Cantrelletal-03,NiWM-11}.

In this paper,  we  consider a $2\times 2$ Lotka-Volterra competition-diffusion system of the form
\begin{equation}
\left\{
\begin{aligned}[cll]
&U_{t} =d_{1}\Delta U+U\left(r_{1}(x)-b_{1}U- c_{1}V\right)\  \ \ \  \ &&\mbox{in}\ \Omega\times\mathbb{R}_{+},\\
&V_{t} =d_{2}\Delta V+V\left(r_{2}(x)-c_{2} V\right)\  \ \ \  \ && \mbox{in}\ \Omega\times\mathbb{R}_{+},\\
&\frac{\partial U}{\partial n}=\frac{\partial V}{\partial n}=0 \  \ \ \  \ & &  \mbox{on}\ \partial\Omega\times\mathbb{R}_{+},\\
&U(0,x)=U_{0}(x),\ \ V(0,x)=V_{0}(x) \ \ \ \ \ \ \ \ \ & &\mbox{in}\ \Omega. \\
\end{aligned}
\right.
\label{control equation}
\end{equation}
where $U(t,x)$ and $V(t,x)$ respectively represent the population densities of
the prey and the predator at location $x$ and time $t>0$,
all model parameters are positive,
the parameters $d_{1}$ and $d_{2}$ are the dispersal rates,
the growth rates $r_{1}(x)$ and $r_{2}(x)$ of the prey and the predator are spatially heterogeneous,
the Laplacian operator $\Delta$ is defined by $\Delta=\sum_{i=1}^N\frac{\partial^{2}}{\partial x_{i}^{2}}$,
the vector  $n$ denotes the outward unit normal vector on $\partial \Omega$,
the habitat $\Omega$ is a bounded region in $\mathbb{R}^{N}$ with smooth boundary $\partial \Omega$,
and $\mathbb{R}_{+}=(0,+\infty)$.
This model indicates that the evolution of the density  $V$ of the predator  is independent of the density $U$ of the prey.

Taking into consideration that  real environments are highly heterogeneous in the abiotic factors,
we assume that  the growth rates $r_{1}(x)$ and $r_{2}(x)$ of the prey and the predator depend on location $x$.
This implies the spatial heterogeneity. More precisely,
throughout this paper
we assume  that the growth rates $r_{i}(x)$ satisfy the following hypothesis:
\begin{enumerate}
\item[{\bf (H1)}]
The growth rates $r_{i}(x)$ are non-negative functions in $C^{\alpha}( \overline{\Omega})$ for $\alpha \in(0,1)$,
$(r_{1}(x),r_{2}(x))$ is not a constant vector in $\Omega$,
and $r_{i}$ satisfy $\overline{r}_{i}>0$, where $\overline{r}_{i}$, $i=1,2$, are defined by
$$\overline{r}_{i}=\frac{1}{|\Omega|}\int_{\Omega} r_{i}(x)\,dx,$$
where  $|\Omega|$ denotes the measure of the bounded region $\Omega$.
\end{enumerate}
No confusion should arise,
when we use a function with a bar without any qualification we mean the average of this function on $\Omega$.

Our goal of this paper is to investigate the global dynamics of the competition-diffusion system  (\ref{control equation})
and then apply the obtained results to a modified Leslie-Gower model.
More precisely, we obtain that  system (\ref{control equation}) always has the trivial steady states $(0,0)$,
and precisely two semi-trivial steady states $(\theta_{d_{1},r_{1}}/b_{1},0)$ and $(0,\theta_{d_{2},r_{2}}/c_{2})$.
By  employing  the principal spectral theory,
we further prove that the steady states $(0,0)$ and  $(\theta_{d_{1},r_{1}}/b_{1},0)$ are always linearly unstable,
that is, the associated  principal eigenvalues are always negative (see, for instance, \cite{NiWM-11}),
and the linear stability of the semi-trivial steady state $(0,\theta_{d_{2},r_{2}}/c_{2})$
changes as the parameters vary in different regions.
A steady state is called to be linearly stable if the associated principal eigenvalue is positive.
If the associated principal eigenvalues is equal to zero,
then it is referred to as the degenerate case. Otherwise, it is called the non-degenerate case.
In the recent work \cite{He-Ni-16},
He and Ni  obtained a powerful result that
the global dynamics of  the classical  $2\times 2$ Lotka-Volterra competition-diffusion system could be determined  by its local dynamics.
Stimulated by the conclusions obtained by He and Ni \cite{He-Ni-16},
an interesting question arises:
\begin{itemize}
\item{\it Is the semi-trivial steady state $(0,\theta_{d_{2},r_{2}}/c_{2})$ of the the competition-diffusion system  (\ref{control equation})
 globally asymptotically stable when it is not linearly unstable\,?}
\end{itemize}
The answer to this question is affirmative.
The proof is mainly based on the  principal spectral theory and the dynamics of the classical single-species logistic model.

As an application,
we employ the results on the dynamics of the competition-diffusion system  (\ref{control equation})
to a $2\times 2$ diffusive predator-prey system
with modified Leslie-Gower and Holling type II functional response
in heterogeneous environment
\begin{equation}
\label{2D-PP-model}
\left\{
\begin{aligned}[cll]
&U_{t} =d_{1}\Delta U+U(r_{1}(x)-b_{1}U)-VP(U)\  \ \ \  \ &&\mbox{in}\ \Omega\times\mathbb{R}_{+},\\
&V_{t} =d_{2}\Delta V+V\left(r_{2}(x)-\frac{a_{2}V}{U+k_{2}}\right)\  \ \ \  \ && \mbox{in}\ \Omega\times\mathbb{R}_{+},\\
&\frac{\partial U}{\partial n}=\frac{\partial V}{\partial n}=0 \  \ \ \  \ & &  \mbox{on}\ \partial\Omega\times\mathbb{R}_{+},\\
&U(0,x)=U_{0}(x),\ \ V(0,x)=V_{0}(x) \ \ \ \ \ \ \ \ \ & &\mbox{in}\ \Omega, \\
\end{aligned}
\right.
\end{equation}
with a Holling type II functional response
\begin{eqnarray*}
P(U)=\frac{a_{1}U}{U+k_{1}},
\end{eqnarray*}
where the predator's numerical response is  modified Leslie-Gower form proposed by Aziz-Alaoui and Daher Okiye in \cite{Aziz-Alaoui-03},
which extends the classical Leslie form  originated by Leslie in \cite{Leslie-48}.
Here all model parameters are positive.

As one of important classes of Lotka-Volterra models,
the predator-prey systems with modified Leslie-Gower and Holling type II functional response were widely investigated by many authors.
For example,
\cite{Aziz-Alaoui-03} studied the global stability of a unique positive equilibrium by the method of Lyapunov functions
and \cite{Zhu-Wang-11} investigated the periodic solutions
for these systems of ordinary differential equations,
\cite{Nindjin-etal-06} obtained a unique globally asymptotically stable positive equilibrium under some conditions for this system with delays,
\cite{Abid-etal-15,Camara-08,Daher-04,Huang-03} investigated
the dynamics of the corresponding reaction-diffusion models in homogeneous environment.
Here we assume that  the growth rates $r_{1}(x)$ and $r_{2}(x)$ of the prey and the predator satisfy the hypothesis {\bf (H1)}.
This implies the spatial heterogeneity.

We will see that system (\ref{2D-PP-model}) has a trivial steady state $(0,0)$,
and two semi-trivial steady states $(\theta_{d_{1},r_{1}}/b_{1},0)$ and $(0,k_{2}\theta_{d_{2},r_{2}}/a_{2})$,
where  $\theta_{d_{i},r_{i}}$ denotes the unique positive steady state
of system (\ref{eq-logistic}) with $d=d_{i}$ and $h=r_{i}$ for each $i=1,2$.
As we know, it is highly difficult to obtain the monotonicity of the semiflow for a diffusive Lotka-Volterra predator-prey system
with a Holling-type functional response in heterogeneous environment.
This causes a big obstacle to analyze the dynamics of Lotka-Volterra predator-prey systems.
To overcome this obstacle,
we apply the obtained results for the competition-diffusion system (\ref{control equation}) together with the method of upper and lower solutions to prove the global stability of steady states
and the persistence of the two species under some suitable conditions.

This paper is organized as follows.
In section 2 we first introduce some results on the global dynamics of the classical single-species logistic model
with diffusion in spatially heterogeneous environment.
In section 3 we study the global dynamics of the competition-diffusion system (\ref{control equation}).
In section 4 we applied the obtained results to a modified Leslie-Gower model.
Some concluding remarks are given in the final section.

\section{Preliminaries}
\label{sec-prelim}
\setcounter{equation}{0}
\setcounter{lemma}{0}
\setcounter{theorem}{0}
\setcounter{remark}{0}

In this section
we introduce some results on the eigenvalues problems associated with the following
single-species logistic model
\begin{equation}
\left\{
\begin{aligned}
&U_{t}=d\Delta U+U(h(x)-U) \ && \mbox{in}\ \Omega\times\mathbb{R}_{+},\\
&\frac{\partial U}{\partial n}=0 \ &&  \mbox{on}\ \partial\Omega\times\mathbb{R}_{+},
\end{aligned}
\right.
\label{eq-logistic}
\end{equation}
where $d>0$ and the function $h:\Omega \to \mathbb{R}$ satisfies the following hypothesis:
\begin{enumerate}
\item[{\bf (H2)}]
the  function $h$ is non-constant, bounded and measurable.
\end{enumerate}
Under the hypothesis {\bf (H2)} the single-species logistic model (\ref{eq-logistic})
has at most a unique positive steady state denoted by $\theta_{d,h}$
(see, for instance, \cite{Lou-06,NiWM-11} or Lemma \ref{lemma-logistic-equation}).
Then this steady state $\theta_{d,h}$ satisfies the following equation
\begin{eqnarray*}
d\triangle \theta_{d,h}+\theta_{d,h}(h(x)-\theta_{d,h})=0.
\end{eqnarray*}
Let both sides of this equation be divided by $\theta_{d,h}$.
And then integrating  over the bounded region $\Omega$, we get
\begin{eqnarray*}
d\int_{\Omega}|\frac{\nabla\theta_{d,h}(x)}{\theta_{d,h}(x)}|^{2}\,dx+\int_{\Omega}(h(x)-\theta_{d,h}(x))\,dx=0.
\end{eqnarray*}
Since the function $h$ is non-constant,
then the steady state $\theta_{d,h}$ is also non-constant,
which yields that for each $d>0$,
\begin{eqnarray}
\label{ineq-h-theta}
\int_{\Omega}h(x)\,dx<\int_{\Omega}\theta_{d,h}\,dx.
\end{eqnarray}
To analyze the stability of the steady state $\theta_{d,h}$ for system (\ref{eq-logistic}),
it is useful to study the eigenvalue problem with indefinite weight:
\begin{equation}
\label{eig-prob-1}
\left\{
\begin{aligned}
&\Delta \phi+\lambda h(x)\phi =0 \ \ \  && \mbox{ in }\ \Omega,\\
&\frac{\partial \phi}{\partial n}=0 \ \ \  && \mbox{ on }\ \partial\Omega.\\
\end{aligned}
\right.
\end{equation}
A constant $\lambda$ is referred  to as a {\it principal eigenvalue} of the eigenvalue problem (\ref{eig-prob-1})
if the problem (\ref{eig-prob-1}) has a positive solution associated with $\lambda$.
Clearly, $\lambda=0$ is always a principal eigenvalue for each function $h$.
The properties of the principal eigenvalues for (\ref{eig-prob-1}) are stated in the next lemma.

\begin{lemma}
{\rm \cite[Theorem 4.2, p.67]{NiWM-11}}
\label{lemma-eigen-problem-1}
Assume that the function $h$ satisfies the hypothesis {\bf (H2)} and changes sign in $\Omega$.
Then the eigenvalue problem (\ref{eig-prob-1}) has a nonzero principal eigenvalue $\lambda_{1}=\lambda_{1}(h)$
if and only if $\int_{\Omega}h(x)\,dx\neq0$.
More precisely, the following statements hold:
\begin{enumerate}
\item[{\bf (i)}]
if $\int_{\Omega}h(x)\,dx>0$, then $\lambda_{1}(h)<0$.

\item[{\bf (ii)}]
if $\int_{\Omega}h(x)\,dx=0$, then $\lambda_{1}(h)=0$ is a unique principle eigenvalue.

\item[{\bf (iii)}]
if $\int_{\Omega}h(x)\,dx<0$, then $\lambda_{1}(h)>0$. Moreover, $\lambda_{1}(h)$ is given by
\begin{eqnarray}
\label{df-lambda-1}
\lambda_{1}(h)=\inf\left\{\frac{\int_{\Omega}|\nabla\phi|^2}{\int_{\Omega}h\phi^{2}}:
         \phi\in H^{1}(\Omega)\ and \int_{\Omega}h\phi^2>0\right\},
\end{eqnarray}
which satisfies the following properties:
\item[{\bf (iv)}]
$\lambda_{1}(h)>\lambda_{1}(k)$ if $h\leq k$ and $h\not\equiv k$ in $\Omega$.

\item[{\bf (v)}]
$\lambda_{1}(h_{m})\to \lambda_{1}(h)$ as $\|h_{m}-h\|_{\infty}\to 0$,
where $\|\cdot\|_{\infty}$ is the essential supremum norm of $L^{\infty}(\Omega)$.
\end{enumerate}

\end{lemma}
The eigenvalue problem (\ref{eig-prob-1}) is closely related to the next eigenvalue problem
\begin{equation}
\label{eig-prob-2}
\left\{
\begin{aligned}
&d\Delta \phi+h(x)\phi+\mu \phi=0 \ \ \  && \mbox{ in }\ \Omega,\\
&\frac{\partial \phi}{\partial n}=0 \ \ \  && \mbox{ on }\ \partial\Omega.\\
\end{aligned}
\right.
\end{equation}
We recall that  the first eigenvalue  $\mu_{1}(d,h)$ of the problem (\ref{eig-prob-2}) is given by
the following variational characterization (see \cite[Formula (4.9), p.69]{NiWM-11})
\begin{eqnarray*}
\mu_{1}(d,h)=\inf\left\{\int_{\Omega}\left(d|\nabla \phi|^2-h\phi^{2}\right)\,dx:
      \int_{\Omega} \phi^2\,dx=1 \mbox{ for } \phi\in H^{1}(\Omega)\right\}.
\end{eqnarray*}
To study  the stability of the steady states for system (\ref{2D-PP-model}),
we need to get the sign of $\mu_{1}(d,h)$.

\begin{lemma}
{\rm \cite[Proposition 4.4, p.69]{NiWM-11}}
\label{lm-eigen-problem-2}
Assume that the function $h$ satisfies the hypothesis {\bf (H2)} and changes sign in $\Omega$.
Let $\mu_{1}(d,h)$ denote the first eigenvalue of the eigenvalue problem (\ref{eig-prob-2}).
Then the following statements hold:
\begin{enumerate}
\item[{\bf (i)}]
if $\int_{\Omega}h(x)\,dx \geq0$ and $h\not \equiv 0$, then $\mu_{1}(d,h)<0$ for each $d>0$.

\item[{\bf (ii)}]
if $\int_{\Omega}h(x)\,dx<0$, then the sign of $\mu_{1}(d,h)$ has the following trichotomies:
\vskip 0.15cm
\begin{enumerate}
\item[{\bf (B1)}]
$\mu_{1}(d,h)<0$ for $d<1/\lambda_{1}(h)$;
\item[{\bf (B2)}]
$\mu_{1}(d,h)=0$ for $d=1/\lambda_{1}(h)$;
\item[{\bf (B3)}]
$\mu_{1}(d,h)>0$ for $d>1/\lambda_{1}(h)$,
where $\lambda_{1}(h)$ is defined by (\ref{df-lambda-1}).
\end{enumerate}
\item[{\bf (iii)}]
the first eigenvalue $\mu_{1}(d,h)$ is strictly increasing and concave in $d$. Furthermore,
\begin{eqnarray*}
\lim_{d\to0+}\mu_{1}(d,h)=\min_{\Omega}(-h), \ \ \
\lim_{d\to+\infty}\mu_{1}(d,h)=\overline{h}=-\frac{1}{|\Omega|}\int_{\Omega} h(x)\,dx.
\end{eqnarray*}

\item[{\bf (iv)}]
if $h\geq k$ and $h\not\equiv k$ in $\Omega$, then $\mu_{1}(d,h)<\mu_{1}(d,k)$.
In particularly, if $h\leq0$ and $h\not\equiv 0$, then $\mu_{1}(d,h)>0$.
\end{enumerate}
\end{lemma}

The results on the global dynamics of system (\ref{eq-logistic}) are summarized in the next lemma.

\begin{lemma}
{\rm \cite{Hutson-etal-95,NiWM-11}}
\label{lemma-logistic-equation}
Assume that the function $h$ satisfies the hypothesis {\bf (H2)}
in the single-species logistic model (\ref{eq-logistic}).
Then the following statements hold:
\begin{enumerate}
\item[{\bf (i)}]
if $\int_{\Omega}h(x)\,dx\geq0$, then for each $d>0$
system (\ref{eq-logistic}) has a unique positive steady state $\theta_{d,h}$,
which is globally asymptotically stable. Moreover,
this steady state $\theta_{d,h}$ satisfies the limits:
\begin{eqnarray}\label{lim-d-0-infty}
\lim_{d\to0+}\theta_{d,h}(x)=h^{+}(x):=\max\left\{h(x),0\right\},\ \ \ \ \   \lim_{d\to+\infty}\theta_{d,h}=\overline{h}.
\end{eqnarray}

\item[{\bf (ii)}]
if $\int_{\Omega}h(x)\,dx<0$ and the function $h$ changes sign in $\Omega$,
then  system (\ref{eq-logistic}) has a unique positive steady state $\theta_{d,h}$ if and only if $ 0<d<1/\lambda_{1}(h)$,
where the constant $\lambda_{1}(h)$ is defined by (\ref{df-lambda-1}).
Moreover, the steady state $\theta_{d,h}$ is globally asymptotically stable
and $\theta_{d,h}(x)\to h^{+}(x)$ as $d\to 0+$.
If $d\geq1/\lambda_{1}(h)$,
then the trivial steady state $0$ is a global attractor of system (\ref{eq-logistic})
in  $\{U\in \mathbb{R}: U\geq 0\}$.

\item[{\bf (iii)}]
\label{itm-3-lemma-logis-eq}
if the function $h$ satisfies $h(x)\leq0$ for each $x\in \Omega$,
then the trivial steady state $0$  is a global attractor of system (\ref{eq-logistic}) in $\mathbb{R}$.

\end{enumerate}
\end{lemma}
By applying the method of upper and lower solutions (see, for instance, \cite{Pao-92,Sattinger-72}),
the existence and uniqueness of the steady state $\theta_{d,h}$ can be established.
The results on the limit behavior of the steady state $\theta_{d,h}$ are obtained by \cite[Lemmas 2.4 and 2.5]{Hutson-etal-95}.
An outline of the proof for this lemma is given in \cite[Section 4.1]{NiWM-11}.

\section{Global dynamics of the competition-diffusion system}
\label{subsec-control}
\setcounter{equation}{0}
\setcounter{lemma}{0}
\setcounter{theorem}{0}
\setcounter{remark}{0}

In this section,
we give the detailed study of the competition-diffusion system (\ref{control equation}).
Assume that the competition-diffusion system  (\ref{control equation})
has a non-negative steady state $(u,v)$ with $u\geq 0$ and $v\geq 0$.
We linearize the corresponding elliptic system of (\ref{control equation}) at  $(u,v)$
and obtain
\begin{equation}
\label{con-eq-linear}
\left\{
\begin{aligned}
&d_{1}\Delta\Phi+\left(r_1(x)-2b_{1}u-c_{1} v \right)\Phi-c_{1}u\Psi+\mu\Phi=0\ && \mbox{ in }\ \Omega,\\
&d_{2}\Delta\Psi+\left(r_2(x)-2c_{2}v \right)\Psi+\mu\Psi=0 && \mbox{ in }\ \Omega,\\
&\frac{\partial \Phi}{\partial n}=\frac{\partial \Psi}{\partial n}=0 \ && \mbox{ on }\ \partial\Omega.\\
\end{aligned}
\right.
\end{equation}
By the {\it Krein-Rutman Theorem} \cite[Theorem 4.2, p.20]{Smith-95},
the linearized system (\ref{con-eq-linear}) has a principal eigenvalue $\mu_1$,
which is simple and has the least real part among all eigenvalue.
The local dynamics of system (\ref{control equation}) at this steady state $(u,v)$
are determined by the principal eigenvalue $\mu_1\in \mathbb{R}$.
If the associated principal eigenvalues $\mu_1$ is equal to zero,
then it is referred to as the degenerate case. Otherwise, it is called the non-degenerate case.
The steady state $(u,v)$ is called {\it linearly stable} (resp. {\it linearly unstable}) if $\mu_{1}$ is positive (resp. negative)
(see, for instance, \cite{NiWM-11}).

Before studying the linear stability of the semi-trivial steady states,
we first make some preparations.
Assume that the rate $d$ and the function $h$ in the single-species logistic model (\ref{eq-logistic})
are in the form $d=d_{2}$ and $h=r_{2}$,
where $r_{2}$ satisfy the conditions in  {\bf (H1)}.
Then by Lemma \ref{lemma-logistic-equation}
the single-species logistic model (\ref{eq-logistic}) has a unique positive steady state $\theta_{d_2,r_2}$.
Let the constants $\alpha$ and $\beta$ be respectively defined by
\begin{eqnarray}
\label{df-S-S}
\alpha= \inf_{d_2>0}\frac{\overline{r}_1}{\overline{\theta}_{d_2,r_2} }, \ \ \ \ \
\beta= \sup\limits_{d_2>0}\sup\limits_{\overline{\Omega}}\frac{r_1(x)}{\theta_{d_2,r_2}(x)}.
\end{eqnarray}
Then we have the following.

\begin{lemma}
\label{lm-S-S}
Assume that the functions $r_{1}$ and $r_{2}$ satisfy the hypothesis {\bf (H1)}.
Let $\theta_{d_2,r_2}$ denote the unique positive steady state of
the single-species logistic model (\ref{eq-logistic}) with $d=d_{2}$ and $h=r_{2}$.
Then the constants $\alpha$ and $\beta$ in (\ref{df-S-S}) satisfy $\alpha<\beta$.
\end{lemma}
{\bf Proof.}
If $r_{2}$ is a constant,
then $r_{2}(x)\equiv \overline{r}_{2}>0$ for each $x\in \Omega$
and $\theta_{r_{2},d_{2}}(x)\equiv \overline{r}_{2}$ for all $x\in \Omega$ and $d_{2}>0$.
By the hypothesis {\bf (H1)}, the function $r_{1}$ is not a constant.
Hence, we have the following:
\begin{eqnarray*}
\beta=\sup\limits_{d_2>0}\sup\limits_{\overline{\Omega}}\frac{r_1(x)}{\theta_{d_2,r_2}(x)}
    =\frac{1}{\overline{\theta}_{d_2,r_2}}\sup\limits_{\overline{\Omega}}r_1(x)
    >\frac{\overline{r}_{1}}{\overline{\theta}_{d_2,r_2} }=\alpha.
\end{eqnarray*}
Thus, the estimate $\alpha<\beta$ holds.

If $r_{2}$ is not a constant,
then $\theta_{d_2,r_2}$ is positive and non-constant in $\Omega$.
Thus we have
\begin{eqnarray}
\label{est-00}
\begin{split}
0&=\int_{\Omega}
\left(r_{1}(x)\overline{\theta}_{d_2,r_2} -\overline{r}_{1}\theta_{d_2,r_2}(x)\right)\,dx\\
&=\int_{\Omega}\overline{\theta}_{d_2,r_2} \theta_{d_2,r_2}(x)
\left(r_{1}(x)/\theta_{d_2,r_2}(x)-\overline{r}_{1}/\overline{\theta}_{d_2,r_2}\right)\,dx,
\end{split}
\end{eqnarray}
which implies
\begin{eqnarray}
\sup_{\overline{\Omega}}\frac{r_{1}(x)}{\theta_{d_2,r_2}(x)}\geq \frac{\overline{r}_{1}}{\overline{\theta}_{d_2,r_2}}.
\label{est-1}
\end{eqnarray}
Otherwise, suppose that  $r_{1}(x)/\theta_{d_2,r_2}(x)<\overline{r}_{1}/\overline{\theta}_{d_2,r_2}$ for $x\in\Omega$.
Then
\begin{eqnarray*}
\int_{\Omega}\overline{\theta}_{d_2,r_2} \theta_{d_2,r_2}(x)
\left(r_{1}(x)/\theta_{d_2,r_2}(x)-\overline{r}_{1}/\overline{\theta}_{d_2,r_2}\right)\,dx<0,
\end{eqnarray*}
which contradicts (\ref{est-00}).
By (\ref{ineq-h-theta}) and (\ref{lim-d-0-infty}),
we obtain that $\overline{\theta}_{\cdot,\,r_2}$ reaches its minimum at infinity
and $\overline{\theta}_{d_{2},\,r_2}$ is not a constant as $d_{2}$ varies. This
together with (\ref{est-1}) yields
\begin{eqnarray*}
\beta=\sup_{d_2>0}\sup_{\overline{\Omega}}\frac{r_1(x)}{\theta_{d_2,r_2}(x)}
     \geq \sup_{d_2>0}\frac{\overline{r}_{1}}{\overline{\theta}_{d_2,r_2}}
     >\inf_{d_2>0}\frac{\overline{r}_{1}}{\overline{\theta}_{d_2,r_2}}
     =\alpha.
\end{eqnarray*}
Therefore, the proof is now complete.
\hfill$\Box$

By Lemma \ref{lemma-logistic-equation},
the second equation in  system (\ref{control equation}) satisfying {\bf (H1)} has a unique positive steady state
$\theta_{d_2,r_2}/c_{2}$.
Let the sets $\widetilde{I}$, $\widetilde{I}_{1}$ and $\widetilde{I}_{2}$ be defined by
\begin{eqnarray}
\widetilde{I}\!\!\!&=&\!\!\!\left\{d_2 \in \mathbb{R}_{+}: \int_{\Omega} \left(r_1(x)-\frac{c_{1}}{c_{2}}\theta_{d_2,r_2}(x)\right)\,dx<0\right\},\label{df-tld-I}\\
\widetilde{I}_{1}\!\!\!&=&\!\!\!\left\{d_2 \in \mathbb{R}_{+}: r_1-\frac{c_{1}}{c_{2}}\theta_{d_2,r_2}\leq0
       \ \mbox{ and }\ r_1-\frac{c_{1}}{c_{2}}\theta_{d_2,r_2}\not\equiv 0 \right\},\label{df-tld-I-1}\\
\widetilde{I}_{2}\!\!\!&=&\!\!\!\left\{d_2\in \widetilde{I} :\ \sup_{\overline{\Omega}}\left(r_1(x)-\frac{c_{1}}{c_{2}}\theta_{d_2,r_2}(x)\right)>0 \right\},
\label{df-tld-I-2}
\end{eqnarray}
respectively.
Then $\widetilde{I}=\widetilde{I}_1\cup \widetilde{I}_2$.
The results on the linear stabilities of the steady states of 
the competition-diffusion system  (\ref{control equation})  are summarized as follows.

\begin{lemma}
\label{thm-stab-contrl}
Assume that the competition-diffusion system  (\ref{control equation}) satisfies the hypothesis {\bf (H1)}.
Then system  (\ref{control equation}) has a trivial steady state $(0,0)$,
and two semi-trivial steady states $(\theta_{d_{1},r_{1}}/b_{1},0)$ and $(0,\theta_{d_{2},r_{2}}/c_{2})$,
where $\theta_{d_{i},r_{i}}$ denotes the unique positive steady state
of system (\ref{eq-logistic}) with $d=d_{i}$ and $h=r_{i}$ for each $i=1,2$.
Furthermore, the following statements hold:
\begin{enumerate}
\item[{\bf (i)}]
\label{itm-1-thm-stab-ctr1}
the trivial steady state $(0,0)$ and the semi-trivial steady state $(\theta_{d_{1},r_{1}}/b_{1},0)$
are both linearly unstable.

\item[{\bf (ii)}]
\label{itm-2-thm-stab-ctr1}
the semi-trivial steady state $(0,\theta_{d_{2},r_{2}}/c_{2})$ is linearly stable for
$(d_{1},d_{2})$ in $\widetilde{\mathcal{D}}_{+}$
and  linearly unstable for $(d_{1},d_{2})$ in $\widetilde{\mathcal{D}}_{-}$,
where the sets $\widetilde{\mathcal{D}}_{\pm}$ are respectively given by
\begin{eqnarray}
\widetilde{\mathcal{D}}_{+}\!\!\!&=&\!\!\! \left\{(d_{1},d_{2})\in\mathbb{R}_{+}^{2}:
                          \mu_{1}\left(d_1,r_{1}-\frac{c_{1}}{c_{2}}\theta_{d_{2},r_{2}}\right)>0\right\},\label{df-WD-D+}\\
\widetilde{\mathcal{D}}_{-}\!\!\!&=&\!\!\! \left\{(d_{1},d_{2})\in\mathbb{R}_{+}^{2}:
                           \mu_{1}\left(d_1,r_{1}-\frac{c_{1}}{c_{2}}\theta_{d_{2},r_{2}}\right)<0\right\},\nonumber
\end{eqnarray}

and the set $\widetilde{\mathcal{D}}_{+}$ has the following trichotomies:
\begin{enumerate}
\item[{\bf (T1')}]
if $c_{1}/c_{2}\in (0,\alpha]$,
then $\widetilde{\mathcal{D}}_{+}=\varnothing$;
\item[{\bf (T2')}]
if  $c_{1}/c_{2}\in [\beta,+\infty)$,
then $\widetilde{\mathcal{D}}_{+}=\mathbb{R}_{+}^2$;
\item[{\bf (T3')}]
if $c_{1}/c_{2}\in (\alpha,\beta)$,
then $\widetilde{\mathcal{D}}_{+}=\{(d_1,d_2)\in\mathbb{R}_{+}^2: d_2\in \widetilde{I}, \ d_1>\widetilde{\varphi}(d_2)\}$,
where the function $\widetilde{\varphi}$ is defined by
\begin{equation}
\widetilde{\varphi}(d_2):=
\left\{
\begin{aligned}
&0  & \ \mbox{ for } \ d_2\in \widetilde{I}_1,\\
&(\lambda_1\left({r_{1}-c_{1}\theta_{d_2,r_2}/c_{2}}\right))^{-1} & \  \mbox{ for } \ d_2\in \widetilde{I}_2.
\end{aligned}
\right.
\label{df-varphi}
\end{equation}
\end{enumerate}
\end{enumerate}
\end{lemma}
{\bf Proof.}
It is clear that system (\ref{control equation}) has the trivial steady state $(0,0)$.
By Lemma \ref{lemma-logistic-equation} we obtain that
system (\ref{control equation}) has two semi-trivial steady states
$(\theta_{d_{1},r_{1}}/b_{1},0)$ and $(0,\theta_{d_{2},r_{2}}/c_{2})$.
Thus, the first statement holds.

Consider the linearized system (\ref{con-eq-linear}) with $(u,v)=(0,0)$, that is,
\begin{equation*}
\left\{
\begin{aligned}
&d_{1}\Delta\Phi+r_1(x)\Phi+\mu\Phi=0\ && \mbox{ in }\ \Omega,\\
&d_{2}\Delta\Psi+r_2(x)\Psi+\mu\Psi=0 && \mbox{ in }\ \Omega,\\
&\frac{\partial \Phi}{\partial n}=\frac{\partial \Psi}{\partial n}=0 \ && \mbox{ on }\ \partial\Omega.\\
\end{aligned}
\right.
\end{equation*}
Then by Lemma \ref{lm-eigen-problem-2} {\bf (i)},
we see that $\mu_1(d_{i},r_{i})<0$ for $d_{i}>0$.
This proves that  $(0,0)$ is linearly unstable.

Consider the linearized system (\ref{con-eq-linear}) with $(u,v)=(\theta_{d_{1},r_{1}}/b_{1},0)$, that is,
\begin{equation}
\label{eq-linear-semi-1}
\left\{
\begin{aligned}
&d_{1}\Delta\Phi+\left(r_1(x)-2\theta_{d_{1},r_{1}}(x)\right)\Phi
-\frac{c_1}{b_{1}}\theta_{d_{1},r_{1}}(x)\Psi+\mu\Phi=0\ && \mbox{ in }\ \Omega,\\
&d_{2}\Delta\Psi+r_2(x)\Psi+\mu\Psi=0 && \mbox{ in }\ \Omega,\\
&\frac{\partial \Phi}{\partial n}=\frac{\partial \Psi}{\partial n}=0 \ && \mbox{ on }\ \partial\Omega.\\
\end{aligned}
\right.
\end{equation}
By {\bf (i)} and {\bf (iv)} in Lemma \ref{lm-eigen-problem-2},
we obtain that $\mu_1(d_{2},r_{2})<0$ for each $d_{2}>0$,
and $\mu_{1}(d_{1},r_1-2\theta_{d_{1},r_{1}})>\mu_{1}(d_{1},r_1-\theta_{d_{1},r_{1}})=0$ for each $d_{1}>0$.
Then for $\mu=\mu_1(d_{2},r_{2})<0$,
the operator
\begin{eqnarray*}
d_{1}\Delta\cdot+\left(r_1(x)-2\theta_{d_{1},r_{1}}(x)\right)\cdot+\mu\,\cdot
\end{eqnarray*}
is invertible,
which implies that  the eigenvalue problem (\ref{con-eq-linear}) has a negative eigenvalue $\mu_1(d_{2},r_{2})$.
This yields that $(\theta_{d_{1},r_{1}}/b_{1},0)$ is linearly unstable.
Thus, {\bf (i)} is proved.

Consider the linearized system (\ref{con-eq-linear}) with  $(u,v)=(0,\theta_{d_{2},r_{2}}/c_{2})$,
that is,
\begin{eqnarray}
\left\{
\begin{aligned}
&d_{1}\Delta\Phi+\left(r_1(x)-\frac{c_{1}}{c_{2}}\theta_{d_{2},r_{2}}(x)\right)\Phi+\mu\Phi=0\ && \mbox{ in }\ \Omega,\\
&d_{2}\Delta\Psi+(r_2(x)-2\theta_{d_{2},r_{2}}(x))\Psi+\mu\Psi=0 && \mbox{ in }\ \Omega,\\
&\frac{\partial \Phi}{\partial n}=\frac{\partial \Psi}{\partial n}=0 \ && \mbox{ on }\ \partial\Omega.\\
\end{aligned}
\right.
\label{eq-linear-semi-2}
\end{eqnarray}
By similar method used in the proof for the stability of the semi-trivial steady state $(\theta_{d_{1},r_{1}}/b_{1},0)$,
the eigenvalue problem (\ref{eq-linear-semi-2}) with $(d_{1},d_{2})\in\widetilde{\mathcal{D}}_{-}$  has a negative eigenvalue
$$\mu=\mu_{1}\left(d_1,r_{1}-\frac{c_{1}}{c_{2}}\theta_{d_{2},r_{2}}\right)<0.$$
Thus, $(0,\theta_{d_{2},r_{2}}/c_{2})$ is unstable
for $(d_{1},d_{2})\in\widetilde{\mathcal{D}}_{-}$.

We next consider the case that $(d_{1},d_{2})$ is in $\widetilde{\mathcal{D}}_{+}$.
Let $(\Phi_{0},\Psi_{0}) \neq (0,0)$ be a solution of the eigenvalue problem (\ref{eq-linear-semi-2}) \
with some eigenvalue $\mu$.
If  $\Phi_{0}$ satisfies  $\Phi_{0}\neq 0$,
then
\begin{eqnarray*}
\mu\geq \mu_{1}\left(d_1,r_{1}-\frac{c_{1}}{c_{2}}\theta_{d_{2},r_{2}}\right)>0.
\end{eqnarray*}
If  $\Phi_{0}=0$ and $\Psi_{0}\neq 0$,
then
$$\mu\geq \mu_{1}(d_{2},r_2-2\theta_{d_{2},r_{2}})>\mu_{1}(d_{2},r_2-\theta_{d_{2},r_{2}})=0.$$
Thus, $(0,\theta_{d_{2},r_{2}}/c_{2})$ is linearly stable for $(d_{1},d_{2})\in\widetilde{\mathcal{D}}_{+}$.

To obtain the explicit expression of the set $\widetilde{\mathcal{D}}_{+}$,
we first recall that the constants $\alpha$ and $\beta$ defined by (\ref{df-S-S}) satisfies $\alpha<\beta$.
If $c_{1}/c_{2}\in (0,\alpha]$,
then by (\ref{df-S-S}) we have that  for each $d_{2}>0$,
\begin{eqnarray*}
\int_{\Omega}\left(r_{1}(x)-\frac{c_1}{c_2}\theta_{d_2,r_2}(x)\right)\,dx\geq0.
\end{eqnarray*}
Lemma \ref{lm-eigen-problem-2} {\bf (i)} yields that
$\mu_1\left(d_1, r_{1}-\frac{c_1}{c_2}\theta_{d_2,r_2}\right)\leq0$ for all $(d_1,d_2)\in\mathbb{R}_{+}^2$.
Then {\bf (T1')} is proved.

If $c_1/c_2\in [\beta,+\infty)$,
then from (\ref{df-S-S}) it follows that for each $x\in \overline{\Omega}$ and each $d_{2}>0$,
\begin{eqnarray}
r_{1}(x)-\frac{c_1}{c_2}\theta_{d_2,r_2}(x) \leq 0.
\label{est-2}
\end{eqnarray}
We further claim that if $c_1/c_2\in [\beta,+\infty)$, then for all  $d_{2}>0$,
\begin{eqnarray}
r_{1}-\frac{c_1}{c_2}\theta_{d_2,r_2} \not\equiv 0.
\label{est-3}
\end{eqnarray}
In fact, in case $c_1/c_2\in (\beta,+\infty)$, it is clear that the claim holds.
In case $c_1/c_2=\beta$, the proof is divided into three different cases:
{\bf (E1)} if  $r_{1}$ is not a constant and $r_{2}$ is a constant,
then $\theta_{d_{2},r_{2}}$ is a constant for all $d_{2}>0$,
and $r_{1}-\frac{c_1}{c_2}\theta_{d_2,r_2}$ is not a constant.
Thus, the claim holds.
{\bf (E2)}
if $r_{1}$ is a constant and $r_{2}$ is not a constant,
then $\theta_{d_{2},r_{2}}$ is not a constant for all $d_{2}>0$,
which yields the claims.
{\bf (E3)}
if both $r_{i}$ are not constants,
then $\theta_{d_{2},r_{2}}$ is not a constant for all $d_{2}>0$.
Suppose that there exists a $\widehat{d}_{2}$ with $0<\widehat{d}_{2}<+\infty$ such that
$r_{1}-\frac{c_1}{c_2}\theta_{\widehat{d}_{2},r_2}\equiv 0$.
Then for each $x\in \Omega$,
\begin{eqnarray*}
\frac{r_1(x)}{\theta_{\widehat{d}_{2},r_2}(x)}=\frac{c_1}{c_2}=\beta
   = \sup\limits_{d_2>0}\sup_{\overline{\Omega}}\frac{r_1(x)}{\theta_{d_2,r_2}(x)},
\end{eqnarray*}
which yields that $\theta_{\widehat{d}_{2},r_2}\leq\theta_{d_2,r_2}$
and $\theta_{\widehat{d}_{2},r_2} \not\equiv \theta_{d_2,r_2}$ for all $d_{2}>0$.
This is a contradiction with the fact that
(\ref{ineq-h-theta}) and (\ref{lim-d-0-infty}) hold if $r_{2}$ is not a constant.
Hence, the claim is proved.
Thus, by (\ref{est-2}), (\ref{est-3}) and Lemma \ref{lm-eigen-problem-2} {\bf (iv)}.
we obtain {\bf (T2')}.

If $c_1/c_2\in (\alpha,\beta)$,
then by (\ref{ineq-h-theta}), (\ref{lim-d-0-infty}) and (\ref{df-S-S}),
there exist $\widetilde{d}_{2}>0$ and $x_{0}\in \overline{\Omega}$ such that
\begin{eqnarray*}
r_{1}(x_{0})-\frac{c_1}{c_2}\theta_{\widetilde{d}_{2},r_{2}}(x_{0})>0
\ \mbox{ and }\ \int_{\Omega}\left(r_1(x)-\frac{c_1}{c_2}\theta_{\widetilde{d}_2,r_2}(x)\right)\,dx<0,
\end{eqnarray*}
which yields that $\widetilde{I}_{2}$ is not an empty set.
By Lemma \ref{lm-eigen-problem-2} {\bf (ii)},
then for either $d_1>0$ and $d_2\in \widetilde{I}_1$ or  $d_1>1/\lambda(r_1-\frac{c_1}{c_2}\theta_{d_2,r_2})$ and $d_2\in \widetilde{I}_2$,
we have $\mu_{1}(d_1,r_{1}-\frac{c_1}{c_2}\theta_{d_{2},r_{2}})>0$.
Thus, {\bf (T3')} holds.
Therefore, the proof is now complete.
\hfill$\Box$

We observe that  the non-trivial steady state $(u,v)$ of the competition-diffusion system (\ref{control equation})
satisfies that $v=\theta_{d_{2},r_{2}}/c_{2}$ and $u$ is the solution of the following single-species logistic model
\begin{equation}
\left\{
\begin{aligned}[cll]
&U_{t} =d_{1}\Delta U+U\left(\widetilde{r}_{1}(x)-b_{1}U\right)\  \ \ \  \ &&\mbox{in}\ \Omega\times\mathbb{R}_{+},\\
&\frac{\partial U}{\partial n}=0 \  \ \ \  \ & &  \mbox{on}\ \partial\Omega\times\mathbb{R}_{+},
\end{aligned}
\right.
\label{eq-u-direct}
\end{equation}
where the function $\widetilde{r}_{1}$ is in the form
$\widetilde{r}_{1}(x)=r_1(x)-c_{1}\theta_{d_2,r_2}(x)/c_{2}$ for $x\in \Omega$.
Let the sets $\widetilde{\mathcal{D}}_{\pm}$ be defined as in Lemma \ref{thm-stab-contrl},
and $\widetilde{\mathcal{D}}_{0}$ are given by
\begin{eqnarray}
\label{df-tild-D-0}
\widetilde{\mathcal{D}}_{0}\!\!\!&=&\!\!\!\left\{(d_{1},d_{2})\in\mathbb{R}_{+}^{2}:
                          \mu_{1}\left(d_1,r_{1}-\frac{c_{1}}{c_{2}}\theta_{d_{2},r_{2}}\right)=0\right\}.
\end{eqnarray}
In order to study the degenerate case,
we next give a description of the set $\widetilde{\mathcal{D}}_{0}$.
\begin{lemma}
\label{lm-D-vo}
Assume that the functions $r_{1}$ and $r_{2}$ satisfy the hypothesis {\bf (H1)}.
Let the set $\widetilde{\mathcal{D}}_{0}$ be defined by (\ref{df-D0})
and a family of sets $\mathcal{E}_{s}$, $s>0$, be in the form
\begin{eqnarray}
\label{df-E-s}
\mathcal{E}_{s}=\{(d_{1},d_2)\in \mathbb{R}^{2}_{+}: r_1\equiv s\theta_{d_2,r_2}\}, \ \  \ s>0.
\end{eqnarray}
Then for each $s>0$,
the set $\mathcal{E}_{s}$
has the following dichotomies:
\begin{enumerate}
\item[{\bf (D1)}]
 $\mathcal{E}_{s}=\varnothing$;
\item[{\bf (D2)}]
$\mathcal{E}_{s}=\{(d_{1},d_{2}^{*})\in \mathbb{R}_{+}^{2}:d_1>0\}$,
where $d_{2}^{*}$ is the only $d_2$ satisfying $r_1\equiv s\theta_{d_2,r_2}$,
\end{enumerate}
and the set $\widetilde{\mathcal{D}}_{0}$ has the following expression
\begin{equation}
\widetilde{\mathcal{D}}_{0}=\left\{
\begin{aligned}
&\,\varnothing &&\ \ \mbox{ if }\ \frac{c_1}{c_2}\in (0,\alpha)\cup[\beta,+\infty),\\
&\left\{(d_1,d_2)\in\mathbb{R}_{+}^{2}: r_1\equiv \alpha\theta_{d_2,r_2}\right\} &&\ \ \mbox{ if }\ \frac{c_{1}}{c_{2}}=\alpha,\\
&\partial\widetilde{\mathcal{D}}_{+}\cup \left\{(d_1,d_2)\in\mathbb{R}_{+}^{2}: r_1\equiv \frac{c_{1}}{c_{2}} \theta_{d_2,r_2}\right\}
   &&\ \ \mbox{ if }\ \frac{c_{1}}{c_{2}} \in(\alpha,\beta),\\
\end{aligned}
\right.
\label{new-sigam0}
\end{equation}
where $\widetilde{\mathcal{D}}_{+}$ is in the form (\ref{df-WD-D+}), the boundary $\partial\widetilde{\mathcal{D}}_{+}$ of $\widetilde{\mathcal{D}}_{+}$ is given by 
\begin{equation}
\partial\widetilde{\mathcal{D}}_{+}=\left\{
\begin{aligned}
&\varnothing &\qquad \mbox{if} \; d_2\in \widetilde{I}_{1},\\
&\left\{(d_1,d_2)\in\mathbb{R}_{+}^{2} :
d_1=\left(\lambda_1\left(r_{1}-\frac{c_1}{c_2}\theta_{d_2,r_2}\right)\right)^{-1}\right\}
&\qquad \mbox{if} \; d_2\in \widetilde{I}_{2}.
\end{aligned}
\right.
\label{sigam0-case}
\end{equation}
and the sets $\widetilde{I}_{1}$ and $\widetilde{I}_{2}$ are defined by (\ref{df-tld-I-1}) and (\ref{df-tld-I-2}), respectively.
\label{lm-desp-sigma-0}
\end{lemma}
{\bf Proof.}
The proof for the dichotomies {\bf (D1)} and {\bf (D2)} is divided into three different cases:
{\bf (E1)} $r_{1}$ is a constant and $r_{2}$ is non-constant;
{\bf (E2)} $r_{1}$ is non-constant and $r_{2}$ is a constant;
{\bf (E3)} both $r_{1}$ and $r_{2}$ are non-constant.

{\bf Case (E1)}:
If $r_{1}$ is a constant and $r_{2}$ is non-constant,
then $\theta_{d_{2},r_{2}}$ is non-constant,
which yields that no constants $d_{2}$ satisfy $r_1\equiv s\theta_{d_2,r_2}$ for each $s>0$.
Then $\mathcal{E}_{s}=\varnothing$.

{\bf Case (E2)}:
If $r_{1}$ is non-constant and $r_{2}$ is a constant,
then $\theta_{d_{2},r_{2}}$ is a constant.
This implies that for each $s>0$,
there are also no constants $d_{2}$ such that $r_1\equiv s\theta_{d_2,r_2}$ holds.
Then $\mathcal{E}_{s}=\varnothing$.

{\bf Case (E3)}:
If both $r_{1}$ and $r_{2}$ are non-constant,
then $\theta_{d_{2},r_{2}}$ is also non-constant
and satisfies
\begin{eqnarray}
d_{2}\Delta \theta_{d_2,r_2}+\theta_{d_2,r_2}(r_{2}(x)-\theta_{d_2,r_2})=0.
\label{eq-d2}
\end{eqnarray}
For each $s>0$,
assume that $(d_{1},d_{2})$ is in the set $\mathcal{E}_{s}$.
Then $d_{2}$ satisfies $r_1\equiv s\theta_{d_2,r_2}$.
This together with (\ref{eq-d2}) yields that
\begin{eqnarray*}
d_{2}=\frac{\int_{\Omega}\left(\theta_{d_{2},r_{2}}(x)-r_{2}(x)\right)\,dx}
      {\int_{\Omega}|(\nabla\theta_{d_{2},r_{2}}(x))/\theta_{d_{2},r_{2}}(x)|^{2}dx}
     =\frac{\int_{\Omega}(r_{1}(x)/s-r_{2}(x))dx}{\int_{\Omega}|(\nabla r_{1}(x))/r_{1}(x)|^{2}dx}.
\end{eqnarray*}
By (\ref{ineq-h-theta}) we have $d_{2}>0$.
Hence, the positive constant $d_{2}$ with the property $r_1\equiv s\theta_{d_2,r_2}$ is well-defined and unique.
This yields that $\mathcal{E}_{s}$ satisfies {\bf (D2)}.
Thus, {\bf (D1)} and {\bf (D2)} are obtained.

Next we prove (\ref{new-sigam0}).
If $c_{1}/c_{2}\in (0,\alpha)$,
then by (\ref{df-S-S})
we have that
\begin{eqnarray*}
\int_{\Omega}\left(r_1(x)-\frac{c_1}{c_2}\theta_{d_2,r_2}(x)\right)dx>0 \ \ \ \mbox{ for each }\ d_{2}>0.
\end{eqnarray*}
This together with Lemma \ref{lm-eigen-problem-2} {\bf (i)}
yields that for each $d_{1}>0$,
\begin{eqnarray*}
\mu_1\left(d_1,r_1-\frac{c_1}{c_2} \theta_{d_2,r_2}\right)<0.
\end{eqnarray*}
Thus, $\widetilde{\mathcal{D}}_{0}=\varnothing$ for $c_{1}/c_{2}\in (0,\alpha)$.

If $c_{1}/c_{2}\in [\beta,+\infty)$,
then Lemma \ref{thm-stab-contrl} yields $\widetilde{\mathcal{D}}_{+}=\mathbb{R}_{+}^2$.
This implies that  $\widetilde{\mathcal{D}}_{0}=\varnothing$ for $c_{1}/c_{2}\in [\beta,+\infty)$.

If $c_{1}/c_{2}=\alpha$,
then (\ref{df-S-S}) yields that
\begin{eqnarray}
\int_{\Omega}\left(r_1(x)-\alpha\theta_{d_2,r_2}(x)\right)dx\geq0, \ \ \ \mbox{ for each }d_{2}>0.
\label{ineq-03}
\end{eqnarray}
Assume that  $d_{2}$ satisfies  $r_{1}\equiv \alpha \theta_{d_2,r_2}$.
Then $\mu_{1}(d_{1},0)=0$ is the first eigenvalue of the Neumann boundary problem (\ref{eig-prob-2}) with $h\equiv 0$.
This implies
\begin{eqnarray}
\left\{(d_1,d_2)\in\mathbb{R}_{+}^{2}: r_1\equiv \alpha \theta_{d_2,r_2}\right\} \subset \widetilde{\mathcal{D}}_{0}.
\label{subset-relation}
\end{eqnarray}
Assume that $d_{2}$  satisfies $r_{1}\not\equiv \alpha \theta_{d_2,r_2}$.
Then by (\ref{ineq-03}) and Lemma \ref{lm-eigen-problem-2} {\bf (i)},
we obtain that $\mu_1(d_1,r_1-\alpha\theta_{d_2,r_2})<0$ for each $d_{1}>0$.
This together with (\ref{subset-relation}) yields that
$\widetilde{\mathcal{D}}_{0}=\left\{(d_1,d_2)\in\mathbb{R}_{+}^{2}: r_1\equiv \alpha\theta_{d_2,r_2}\right\}$
for $c_{1}/c_{2}=\alpha$.

If $c_{1}/c_{2}\in(\alpha,\beta)$,
then Lemma \ref{thm-stab-contrl} {\bf (ii)} we have
$$
\widetilde{\mathcal{D}}_{+}=\left\{(d_1,d_2)\in\mathbb{R}_{+}^2: d_2\in \widetilde{I}, \ d_1>\widetilde{\varphi}(d_2)\right\}.
$$
Thus $\partial \widetilde{\mathcal{D}}_{+}$  is in the form
\begin{eqnarray*}
\partial \widetilde{\mathcal{D}}_{+}=\left\{(d_1,d_2)\in\mathbb{R}_{+}^2: d_2\in \widetilde{I}, \ d_1=\widetilde{\varphi}(d_2)\right\}.
\end{eqnarray*}
If $d_{2}\in \widetilde{I}_{1}$,
then $\widetilde{\varphi}(d_2)=0$, which yields $\partial \widetilde{\mathcal{D}}_{+}=\varnothing$.
If $d_{2}\in \widetilde{I}_{2}$,
then by  Lemma \ref{lm-eigen-problem-2} {\bf (ii)},
\begin{eqnarray*}
d_{1}=\widetilde{\varphi}(d_2)=\left(\lambda_1\left(r_{1}-\frac{c_1}{c_2}\theta_{d_2,r_2}\right)\right)^{-1},\ \ \
\mu_{1}\left(d_{1},r_{1}-\frac{c_1}{c_2}\theta_{d_2,r_2}\right)=0.
\end{eqnarray*}
Thus, we obtain (\ref{sigam0-case}) and $\partial \widetilde{\mathcal{D}}_{+} \subset \widetilde{\mathcal{D}}_{0}$.

Similarly  to (\ref{subset-relation}),
we obtain
$$
\mathcal{E}_{c_{1}/c_{2}}
=\left\{(d_1,d_2)\in\mathbb{R}_{+}^{2}: r_1\equiv \frac{c_1}{c_2}\theta_{d_2,r_2}\right\} \subset \widetilde{\mathcal{D}}_{0}.
$$
This together with  $\partial \widetilde{\mathcal{D}}_{+} \subset \widetilde{\mathcal{D}}_{0}$ yields
$(\partial\widetilde{\mathcal{D}}_{+}\cup \mathcal{E}_{c_1/c_2})\subset \widetilde{\mathcal{D}}_{0}$ for
$c_1/c_2\in (\alpha,\beta)$.
Thus,
to finish the proof,
it is only necessary to obtain
$\widetilde{\mathcal{D}}_{0}\subset(\partial\widetilde{\mathcal{D}}_{+}\cup \mathcal{E}_{c_1/c_2})$
for $c_1/c_2\in (\alpha,\beta)$.

If $(d_{1},d_{2})\in \widetilde{\mathcal{D}}_{0}$ with $r_{1}\equiv c_{1}\theta_{d_2,r_2}/c_{2}$
then $(d_{1},d_{2})\in \mathcal{E}_{c_1/c_2}$.

If $(d_{1},d_{2})\in \widetilde{\mathcal{D}}_{0}$ with $r_{1}\not\equiv c_1\theta_{d_2,r_2}/c_2$,
then  by Lemma \ref{lm-eigen-problem-2} we have
$$\int_{\Omega} \left(r_{1}-\frac{c_1}{c_2}\theta_{d_2,r_2}\right)dx<0$$
which yields $d_{2}\in \widetilde{I}=\widetilde{I}_1\cup \widetilde{I}_2$.
For $d_{2}\in \widetilde{I}_1$,
Lemma \ref{lm-eigen-problem-2}  yields
$$\mu_{1}\left(d_{1},r_{1}-\frac{c_1}{c_2}\theta_{d_2,r_2}\right)>0.$$
Recall that $\mu_{1}\left(d_1,r_{1}-c_1\theta_{d_{2},r_{2}}/c_2\right)=0$
for $(d_{1},d_{2})\in \widetilde{\mathcal{D}}_{0}$,
then  $d_{2}\in \widetilde{I}_2$.
Following this assertion,
by Lemma \ref{lm-eigen-problem-2} {\bf (ii)}
we have $d_{1}=(\lambda_1({r_{1}-c_1\theta_{d_2,r_2}}/c_2)^{-1}$.
This finishes the proof.
\hfill$\Box$

Finally, we summarize the global dynamics of the competition-diffusion system  (\ref{control equation}) with {\bf (H1)} in the following.
\begin{theorem}
\label{prop-stab-coexist}
Assume that the competition-diffusion system  (\ref{control equation}) satisfies the hypothesis {\bf (H1)}.
Then the following statements hold:
\begin{enumerate}
\item[{\bf (i)}]
\label{itm-1-prp-stab-coexist}
if $(d_{1},d_{2})\in (\widetilde{\mathcal{D}}_{+}\cup \widetilde{\mathcal{D}}_{0})$,
then the semi-trivial steady state $(0,\theta_{d_{2},r_{2}}/c_{2})$ is globally asymptotically stable,
and system (\ref{control equation}) has no positive coexistence steady states.

\item[{\bf (ii)}]
\label{itm-2-prp-stab-coexist}
if $(d_{1},d_{2})\in (\mathbb{R}_{+}^{2}\setminus (\widetilde{\mathcal{D}}_{+}\cup \widetilde{\mathcal{D}}_{0}))$,
then the semi-trivial steady state $(0,\theta_{d_{2},r_{2}}/c_{2})$ is unstable,
and system (\ref{control equation}) has a unique positive coexistence steady state
$(\theta_{d_1,\widetilde{r}_{1}}/b_{1}, \theta_{d_2,r_2}/c_{2})$,
which is globally asymptotically stable.
\end{enumerate}
\end{theorem}
{\bf Proof.}
If $(d_{1},d_{2})\in (\widetilde{\mathcal{D}}_{+}\cup \widetilde{\mathcal{D}}_{0})$,
then according to Lemma \ref{thm-stab-contrl} and  Lemma \ref{lm-D-vo},
we divide the proof for {\bf (\ref{itm-1-prp-stab-coexist})} into three different cases:
$c_{1}/c_{2}=\alpha$, $c_{1}/c_{2}\in (\alpha,\beta)$ and $c_{1}/c_{2}\in [\beta,+\infty)$.
If $(d_{1},d_{2})\in (\mathbb{R}_{+}^{2}\setminus (\widetilde{\mathcal{D}}_{+}\cup \widetilde{\mathcal{D}}_{0}))$,
then by Lemma \ref{thm-stab-contrl} the semi-trivial steady state $(0,\theta_{d_{2},r_{2}}/c_{2})$ is unstable.
The  proof for {\bf (ii)} is divided into the four different cases:
$c_{1}/c_{2}\in (0,\alpha)$, $c_{1}/c_{2}=\alpha$, $c_{1}/c_{2}\in (\alpha,\beta)$ and $c_{1}/c_{2}\in (\beta,+\infty)$.

We only give the proof for the case that
$(d_{1},d_{2})\in (\mathbb{R}_{+}^{2}\setminus (\widetilde{\mathcal{D}}_{+}\cup \widetilde{\mathcal{D}}_{0}))$
and $c_{1}/c_{2}\in (0,\alpha)$,
all other cases can be similarly proved.
By  (\ref{df-S-S}) and $c_{1}/c_{2}\in(0,\alpha)$,
we have that
$\int_{\Omega}\widetilde{r}_{1}(x)\,dx>0$ holds for each $d_{2}>0$.
Hence, applying Lemma \ref{lemma-logistic-equation} to equation (\ref{eq-u-direct}) implies
that system (\ref{control equation}) has a unique positive co-existence steady state
$(\theta_{d_1,\widetilde{r}_{1}}/b_{1}, \theta_{d_2,r_2}/c_{2})$.
To prove that this steady state is globally asymptotically stable,
we first show that it is locally  asymptotically stable.
Substituting $(\theta_{d_1,\widetilde{r}_{1}}/b_{1}, \theta_{d_2,r_2}/c_{2})$ into (\ref{con-eq-linear}) yields
\begin{equation*}
\left\{
\begin{aligned}
&d_{1}\Delta\Phi+\left(\widetilde{r}_{1}(x)-2\theta_{d_1,\widetilde{r}_{1}}(x)\right)\Phi
-\frac{c_{1}}{b_{1}}\theta_{d_1,\widetilde{r}_{1}}(x)\Psi+\mu\Phi=0\ && \mbox{ in }\ \Omega,\\
&d_{2}\Delta\Psi+\left(r_2(x)-2\theta_{d_2,r_2}(x) \right)\Psi+\mu\Psi=0 && \mbox{ in }\ \Omega,\\
&\frac{\partial \Phi}{\partial n}=\frac{\partial \Psi}{\partial n}=0 \ && \mbox{ on }\ \partial\Omega.
\end{aligned}
\right.
\end{equation*}
Then similarly to the proof for Lemma \ref{thm-stab-contrl} {\bf (ii)},
the local stability  is proved.
Choose a sufficiently small $\varepsilon_{0}$
satisfying $\varepsilon_{0}>0$ and $\int_{\Omega}(\widetilde{r}_{1}(x)\pm \varepsilon_{0})\,dx>0$ for each $d_{2}>0$.
Since $v=\theta_{d_2,r_2}/c_{2}$ is a globally asymptotically stable steady state of
the second equation in (\ref{control equation}) with the Neumann boundary condition,
then every solution $(U,V)$ of system (\ref{control equation}) satisfies that
\begin{eqnarray}
\|V(t,\cdot)-\theta_{d_2,r_2}/c_{2}\|_{\infty}\to 0 \ \mbox{ as } \ t\to +\infty.
\label{lim-v-dirct}
\end{eqnarray}
This yields that
for each $\varepsilon$ with $0<\varepsilon<\varepsilon_{0}$,
there exists a sufficiently large $t_{0}$ such that for $t>t_{0}$ and $U\geq 0$,
\begin{eqnarray*}
U\left(\widetilde{r}_{1}(x)-\varepsilon-b_{1}U\right)
\leq U\left(r_{1}(x)-b_{1}U- c_{1}V\right)
\leq U\left(\widetilde{r}_{1}(x)+\varepsilon-b_{1}U\right),
\end{eqnarray*}
and $u=\theta_{d_1,\widetilde{r}_{1}\pm \varepsilon}/b_{1}$ are globally asymptotically steady states
of equation (\ref{eq-u-direct}) with $\widetilde{r}_{1}$ replacing by $\widetilde{r}_{1}\pm \varepsilon$,
then by (2.4) in \cite[p.406]{Lou-06} and the {\it Comparison Principle}
we obtain that for each positive real constant $\delta$
there exists a sufficiently large  $t_{0}$ such that
every solution $(U,V)$ of system (\ref{control equation}) satisfies
$\|U-\theta_{d_1,\widetilde{r}_{1}}/b_{1}\|<\delta$ for each $t>t_{0}$,
together with (\ref{lim-v-dirct}) and the fact that
$(\theta_{d_1,\widetilde{r}_{1}}/b_{1}, \theta_{d_2,r_2}/c_{2})$ is locally asymptotically stable,
yields that $(\theta_{d_1,\widetilde{r}_{1}}/b_{1}, \theta_{d_2,r_2}/c_{2})$ is also globally asymptotically stable.
Therefore, the proof is now complete.
\hfill$\Box$

By this lemma, we have that if the competition-diffusion system  (\ref{control equation}) has a positive coexistence steady state,
then it is globally asymptotically stable.
We also observe that the semi-trivial steady state $(0,\theta_{d_{2},r_{2}}/c_{2})$ is not linearly unstable,
then it is globally asymptotically stable.

\section{Application to a modified Leslie-Gower model}
\label{sec-trivial-state}
\setcounter{equation}{0}
\setcounter{lemma}{0}
\setcounter{theorem}{0}
\setcounter{remark}{0}

In this section we apply the results on the dynamics of the competition-diffusion system  (\ref{control equation})
to the modified Leslie-Gower model (\ref{2D-PP-model}).

Similar to the competition-diffusion system  (\ref{control equation}),
we  define $\mathcal{D}_{\pm}$  and $\mathcal{D}_{0}$ by
\begin{eqnarray}
\mathcal{D}_{+}\!\!\!&:=&\!\!\! \left\{(d_{1},d_{2})\in\mathbb{R}_{+}^{2}:
                          \mu_{1}\left(d_1,r_{1}-\frac{a_{1}k_{2}}{a_{2}k_{1}}\theta_{d_{2},r_{2}}\right)>0\right\},
                          \label{df-D+}\\
\mathcal{D}_{-}\!\!\!&:=&\!\!\! \left\{(d_{1},d_{2})\in\mathbb{R}_{+}^{2}:
                          \mu_{1}\left(d_1,r_{1}-\frac{a_{1}k_{2}}{a_{2}k_{1}}\theta_{d_{2},r_{2}}\right)<0\right\}
                          \nonumber,\\
\mathcal{D}_{0}\!\!\!&:=&\!\!\!\left\{(d_{1},d_{2})\in\mathbb{R}_{+}^{2}:
                          \mu_{1}\left(d_1,r_{1}-\frac{a_{1}k_{2}}{a_{2}k_{1}}\theta_{d_{2},r_{2}}\right)=0\right\}.\label{df-D0}
\end{eqnarray}
The main results on global stability of  the semi-trivial steady state $(0,k_{2}\theta_{d_{2},r_{2}}/a_{2})$
are summarized as follows.

\begin{theorem}
\label{thm-global-1}
Assume that the modified Leslie-Gower model (\ref{2D-PP-model}) satisfies the hypothesis {\bf (H1)},
and the parameters $a_{1}$, $a_{2}$, $k_{1}$ and $k_{2}$ satisfy one of the following conditions:
{\bf (C1)}
$a_1k_2/(a_2k_1)\in[\beta,+\infty)$ and $k_1\geq k_2$;
{\bf (C2)}
$a_1k_2/(a_2k_1)\in(\alpha,\beta)$ and $k_1=k_2$;
{\bf (C3)}
$a_1k_2/(a_2k_1)\in(0,\alpha]$ and $k_{2}\geq k_1$,
where the constants $\alpha$ and $\beta$ are given by (\ref{df-S-S}).
Then  for each $(d_1,d_2)\in \mathcal{D}_{+}\cup \mathcal{D}_{0}$,
the semi-trivial steady state $(0,k_{2}\theta_{d_{2},r_{2}}/a_{2})$ is globally asymptotically stable,
and system (\ref{2D-PP-model}) has no positive steady states.
\end{theorem}

We next state the main results on the the persistence and coexistence states of the modified Leslie-Gower model (\ref{2D-PP-model}).
More precisely,
assume that the modified Leslie-Gower model (\ref{2D-PP-model})  satisfies
one of the following assumptions:
\begin{eqnarray*}
\mbox{{\bf (A1)}}& a_1k_2/(a_2k_1)\in (0,\alpha],\ \ k_1\leq k_2,\ \
(d_1,d_2)\in\mathbb{R}_{+}^2\setminus\left\{(d_1,d_2)\in\mathbb{R}_{+}^2 \ | \ r_1(x)\equiv \alpha\theta_{d_2,r_2}\right\}.\\
\mbox{{\bf (A2)}}& a_1k_2/(a_2k_1)\in (0,\alpha],\ \ k_1=k_2,\ \
(d_1,d_2)\in\mathbb{R}_{+}^2\setminus\left\{(d_1,d_2)\in\mathbb{R}_{+}^2 \ | \ r_1(x)\equiv \alpha\theta_{d_2,r_2}\right\}.\\
\mbox{{\bf (A3)}}&
a_1k_2/(a_2k_1) \in(\alpha,\beta), \ \  k_1= k_2,\ \
(d_{1},d_{2})\in  \left\{(d_1,d_2)\in\mathbb{R}_{+}^2 \ | \ d_2\in I_2, \ d_1<\varphi(d_2)\right\}.
\end{eqnarray*}
where the function $\varphi$ is defined as in Theorem \ref{thm-local}.
By applying the method of upper and lower solutions
and the global dynamics of the competition-diffusion system (\ref{control equation}),
 we prove the existence and the stability of a unique positive steady state, and
establish the uniform persistence of the two species under some suitable conditions.
The modified Leslie-Gower model (\ref{2D-PP-model}) is said to be uniformly persistent \cite[p.61]{Smith-11}
if there exists a positive constant $\delta>0$ such that for each non-negative initial value $(U_{0},V_{0})$
with $U_{0}(x)\not\equiv 0$ and $V_{0}(x)\not\equiv 0$, the solution $(U,V)$ satisfies
\begin{eqnarray*}
\liminf_{t\to+\infty}\min_{x\in \overline{\Omega}} U(t,x)\geq \delta,\ \ \
\liminf_{t\to+\infty}\min_{x\in \overline{\Omega}} V(t,x)\geq \delta.
\end{eqnarray*}
The main results on the the persistence and coexistence states are stated as follows.
\begin{theorem}
\label{thm-coexist}
Assume that that the modified Leslie-Gower model (\ref{2D-PP-model})  satisfies the hypothesis {\bf (H1)}.
Then the following statements hold:
\begin{enumerate}
\item[{\bf (i)}] if either {\bf (A1)} or {\bf (A3)} holds,
then system (\ref{2D-PP-model}) is uniformly persistent.

\item[{\bf (ii)}] if either {\bf (A2)} or {\bf (A3)} holds,
then system (\ref{2D-PP-model}) has at most one positive coexistence steady state.
Further, if there exists a positive coexistence steady state,
then it is a unique positive coexistence steady state, and is globally asymptotically stable.
\end{enumerate}
\end{theorem}

\subsection{Proof of Theorem \ref{thm-global-1}}
To prove Theorem \ref{thm-global-1},
we first make some preparations and start by the explicit expressions of $\mathcal{D}_{+}$ and $\mathcal{D}_{0}$.
Let the sets $I$, $I_{1}$ and $I_{2}$ be denoted by
\begin{eqnarray}
I\!\!\!&:=&\!\!\!\left\{d_2 \in \mathbb{R}_{+}: \int_{\Omega} \,\left(r_1(x)-\frac{a_1k_2}{a_2k_1}\theta_{d_2,r_2}(x)\right)\,dx<0\right\},
          \nonumber\\
I_{1}\!\!\!&:=&\!\!\!\left\{d_2 \in \mathbb{R}_{+}: r_1-\frac{a_1k_2}{a_2k_1}\theta_{d_2,r_2}\leq0
       \ \mbox{ and }\ r_1-\frac{a_1k_2}{a_2k_1}\theta_{d_2,r_2}\not\equiv 0 \right\}, \label{df-I-1}\\
I_{2}\!\!\!&:=&\!\!\!\left\{d_2\in I :\ \sup_{\overline{\Omega}}\left({r_{1}(x)-\frac{a_1k_2}{a_2k_1}\theta_{d_2,r_2}}(x)\right)>0 \right\},
           \label{df-I-2}
\end{eqnarray}
respectively.
Then $I=I_1\cup I_2$.
By the similar method used in Lemma \ref{thm-stab-contrl} {\bf (ii)} we have the following lemma.
\begin{lemma}
\label{thm-local}
Assume that the modified Leslie-Gower model (\ref{2D-PP-model}) satisfies the hypothesis {\bf (H1)}.
Then the set $\mathcal{D}_{+}$ has the following trichotomies:
\begin{enumerate}
\item[{\bf (T1)}]
if $a_1k_2/(a_2k_1)\in (0,\alpha]$, then $\mathcal{D}_{+}=\varnothing$;
\item[{\bf (T2)}]
if  $a_1k_2/(a_2k_1)\in [\beta,+\infty)$, then $\mathcal{D}_{+}=\mathbb{R}_{+}^2$;
\item[{\bf (T3)}]
if $a_1k_2/(a_2k_1)\in (\alpha,\beta)$, then
$\mathcal{D}_{+}=\{(d_1,d_2)\in\mathbb{R}_{+}^2: d_2\in I, \ d_1>\varphi(d_2)\}$,
where the function $\varphi$ is defined by
\begin{equation*}
\varphi(d_2):=
\left\{
\begin{aligned}
&0  & \ \mbox{ for } \ d_2\in I_1,\\
&\left(\lambda_1\left({r_{1}-\frac{a_1k_2}{a_2k_1}\theta_{d_2,r_2}}\right)\right)^{-1} & \  \mbox{ for } \ d_2\in I_2,
\end{aligned}
\right.
\end{equation*}
\end{enumerate}
and $\mathcal{D}_{0}=\widetilde{\mathcal{D}}_{0}$ with $c_{1}/c_{2}$ replaced by $a_{1}k_{2}/(a_{2}k_{1})$,
where $\widetilde{\mathcal{D}}_{0}$ is defined as in Lemma \ref{lm-D-vo}.
\end{lemma}

We next consider the invariant regions of system (\ref{2D-PP-model}).
\begin{lemma}
\label{lm-invar-reg}
Assume that the modified Leslie-Gower model (\ref{2D-PP-model}) satisfies the hypothesis {\bf (H1)}.
Let the constants $M_{1}$ and $M_{2}$ be defined by
$M_1=\sup\limits_{\overline{\Omega}}r_1(x)$
and $M_2=\sup\limits_{\overline{\Omega}}r_2(x)$, respectively.
Then the set $\mathcal{A}$ of the form
\begin{eqnarray*}
\mathcal{A}:=\left\{(U,V): 0\leq U \leq \frac{M_1}{b_1} \; ,\; 0\leq V \leq \frac{(M_1+b_{1}k_2)M_2}{a_2b_{1}} \right\}
\end{eqnarray*}
is an invariant region of system (\ref{2D-PP-model}).
Furthermore, the set  $\mathcal{A}$ is a global attractor of system (\ref{2D-PP-model})
in the set $\{(U,V)\in\mathbb{R}^{2}: U\geq 0,\, V\geq 0\}$.
\end{lemma}
{\bf Proof.}
Similarly to \cite[Lemma 5]{Daher-04},
we obtain that the unique solution of system (\ref{2D-PP-model}) is nonnegative and defined in the set $[0,+\infty)\times \Omega$.
Then by the {\it Comparison Principle},
we obtain that the set $\mathcal{A}$ is an invariant region of system (\ref{2D-PP-model}),
and $U(t,x)\leq\widetilde{U}(t)$ for $t\geq 0$ and $x\in \Omega$,
where $\widetilde{U}$ is the solution of the initial value problem
\begin{eqnarray*}
\frac{d\widetilde{U}}{dt}=\widetilde{U}(M_1-b_1\widetilde{U}), \ \  \
\widetilde{U}(0)=\max\limits_{\overline \Omega}U(x,0).
\end{eqnarray*}
Hence, for $t\geq 0$ and $x\in \Omega$,
\begin{eqnarray}
U(t,x)\leq\widetilde{U}(t)=\frac{\widetilde{U}(0) M_1}{(M_1-b_1\widetilde{U}(0))e^{-M_1t}+b_1\widetilde{U}(0)},
\label{est-U}
\end{eqnarray}
which yields that  for each $U_{0}$ with $0\leq U_{0} \leq M_1/b_1$,
the solution $(U,V)$ of system (\ref{2D-PP-model}) satisfies $0\leq U \leq M_1/b_1$.
Further, consider the following equation
\begin{eqnarray*}
\frac{d\widetilde{V}}{dt}=\widetilde{V}\left(M_2-\frac{a_2b_1}{M_1+k_2b_1}\widetilde{V}\right), \ \ \
\widetilde{V}(0)=\max\limits_{\overline \Omega}V(x,0).
\end{eqnarray*}
Then for each $(U_{0}, V_{0})\in \mathcal{A}$,
we obtain that for each $t\geq 0$ and each $x\in \Omega$,
\begin{eqnarray}
V(t,x)\leq\widetilde{V}(t)=\frac{\widetilde{V}(0)M_{2}}{\left(M_{1}-K_{V}\widetilde{V}(0)\right)e^{-M_{2}t}+K_{V}\widetilde{V}(0)},
\label{est-V}
\end{eqnarray}
where $K_{V}=a_{2}b_{1}/(M_{1}+b_{1}k_{2})$.
The attraction of the set $\mathcal{A}$ is obtained by letting $t\to +\infty$ in (\ref{est-U}) and (\ref{est-V}).
Therefore, the proof is now complete.
\hfill$\Box$

Now we are in the position to prove Theorem \ref{thm-global-1} by applying
the method of upper and lower solutions together with Lemma \ref{prop-stab-coexist}.

{\bf Proof of Theorem \ref{thm-global-1}.}
By Lemma \ref{lm-invar-reg} it suffices to consider the asymptotical behaviors of all solutions of  system (\ref{2D-PP-model})
with $0\leq U_{0}(x)\leq M_{1}/b_{1}$ and $0\leq V_{0}(x)\leq (M_{1}+b_{1}k_{2})M_{2}/(a_{2}b_{1})$
for each $x\in \Omega$.
Thus in the following we always assume that $(U_{0},V_{0})\in \mathcal{A}$.
The detailed proof is divided into three different cases:
{\bf (E1')} $a_1k_2/(a_2k_1)\in[\beta,+\infty)$ and $k_1\geq k_2$;
{\bf (E2')} $a_1k_2/(a_2k_1)\in(\alpha,\beta)$, $k_1=k_2$ and $(d_1,d_2)\in \mathcal{D}_{+}\cup \mathcal{D}_{0}$;
{\bf (E3')} $a_1k_2/(a_2k_1)\in(0,\alpha]$, $k_{2}\geq k_1$ and $(d_1,d_2)\in \mathcal{D}_{+}\cup \mathcal{D}_{0}$.

{\bf Case (E1')}:
If $a_1k_2/(a_2k_1)\in[\beta,+\infty)$ and $k_1\geq k_2$,
then by Theorem \ref{thm-local} we have $\mathcal{D}_{+}=\mathbb{R}_{+}^2$.
Consider the following systems:
\begin{equation}
\label{iterat-1}
\left\{
\begin{aligned}[cll]
&\frac{\partial U_{1j}}{\partial t}
=d_{1}\Delta U_{1j}+U_{1j}\left(r_{1}(x)-b_{1}U_{1j}-\gamma_{1j}V_{1j}\right)
\  \ \ \  \ &&\mbox{in}\ \Omega\times\mathbb{R}_{+},\\
&\frac{\partial V_{1j}}{\partial t}  =d_{2}\Delta V_{1j}+V_{1j}\left(r_{2}(x)-\eta_{1j}V_{1j}\right)
\  \ \ \  \ && \mbox{in}\ \Omega\times\mathbb{R}_{+},\\
&\frac{\partial U_{1j}}{\partial n}=\frac{\partial V_{1j}}{\partial n}=0
\  \ \ \  \ & &  \mbox{on}\ \partial\Omega\times\mathbb{R}_{+},\\
&U_{1j}(0,x)=U_{0}(x),\ \ V_{1j}(0,x)=V_{0}(x) \ \ \ \ \ \ \ \ \ & &\mbox{in}\ \Omega. \\
\end{aligned}
\right.
\end{equation}
where $j=1,2$.
Define the constants $\gamma_{1j}$ and $\eta_{1j}$, $j=1,2$, by
\begin{eqnarray*}
\gamma_{11}=\frac{a_{1}}{k_{1}}, \ \ \ \eta_{11}=\frac{a_{2}}{k_{2}},\ \ \
\gamma_{12}=\frac{a_{1}b_{1}}{M_{1}+b_{1}k_{1}}, \ \ \ \eta_{12}=\frac{a_{2}b_{1}}{M_{1}+b_{1}k_{2}}.
\end{eqnarray*}
then by $a_1k_2/(a_2k_1)\geq \beta$ and $k_2\leq k_1$ we have that
$\gamma_{11}/\eta_{11}\geq \beta$ and $\gamma_{12}/\eta_{12}\geq \beta$.
This together with Lemma \ref{thm-stab-contrl} {\bf (ii)} and Lemma \ref{prop-stab-coexist} {\bf (i)}
yields that for each $j=1$,
system (\ref{iterat-1}) has a globally asymptotically stable steady state $(0,\theta_{d_{2},r_{2}}/\eta_{1j})$
and has no positive coexistence steady states.
Clearly, the solution $(U,V)$ of system (\ref{2D-PP-model})
is an upper solution of system (\ref{iterat-1}) with $j=1$
and a lower solution of system (\ref{iterat-1}) with $j=2$
(see, for instance, \cite[p.21]{Pao-92}).
Then we obtain that
\begin{eqnarray*}
U_{11}(t,x)\leq U(t,x)\leq U_{12}(t,x),\ \
V_{11}(t,x)\leq V(t,x)\leq V_{12}(t,x), \ \ t>0, \ x\in \Omega.
\end{eqnarray*}
Thus for each small $\varepsilon>0$,
there exists a time $t_{1}>0$ such that $0\leq U(t,x)<\varepsilon$
for each $t\geq t_{1}$ and $x\in \Omega$,
which implies that for $t\geq t_{1}$,
\begin{eqnarray*}
V_{11}\left(r_{2}(x)-\eta_{11}V_{11}\right)
\leq V\left(r_{2}(x)-\frac{a_{2}V}{U+k_{2}}\right)
\leq V\left(r_{2}(x)-\frac{a_{2}V}{\varepsilon+k_{2}}\right),\ \ t\geq t_{1}, \ \ x\in \Omega.
\end{eqnarray*}
By applying the method of upper and lower solutions again,
we obtain that
\begin{eqnarray*}
\frac{k_{2}}{a_{2}}\theta_{d_{2},r_{2}}(x) \leq \liminf_{t\to+\infty} V(t,x)
\leq \limsup_{t\to+\infty}V(t,x)\leq \frac{k_{2}+\varepsilon}{a_{2}}\theta_{d_{2},r_{2}}(x), \ \ x\in \Omega.
\end{eqnarray*}
Then we obtain that  $(U(t,x), V(t,x))\to (0,k_{2}\theta_{d_{2},r_{2}}/a_{2})$ as $t\to +\infty$.
Then this case  is proved.

{\bf Case (E2')}: If $a_1k_2/(a_2k_1)\in(\alpha,\beta)$, $k_1=k_2$ and $(d_1,d_2)\in \mathcal{D}_{+}\cup \mathcal{D}_{0}$,
then
\begin{eqnarray*}
\frac{\gamma_{11}}{\eta_{11}}=\frac{\gamma_{12}}{\eta_{12}}=\frac{a_{1}}{a_{2}}\in (\alpha,\beta),
\end{eqnarray*}
which together with Lemmas \ref{lm-D-vo}, \ref{thm-stab-contrl} and \ref{prop-stab-coexist},
yields that for each $j=1$,
system (\ref{iterat-1}) has a globally asymptotically stable steady state $(0,\theta_{d_{2},r_{2}}/\eta_{1j})$
and has no positive coexistence steady states.
Thus similarly to the case {\bf (E1')}, we can prove {\bf (E2')}.

{\bf Case (E3')}: If $a_1k_2/(a_2k_1)\in(0,\alpha]$, $k_{2}\geq k_1$ and $(d_1,d_2)\in \mathcal{D}_{+}\cup \mathcal{D}_{0}$,
then
\begin{eqnarray}
\label{est-4}
0<\frac{\gamma_{11}}{\eta_{11}}=\frac{a_1k_2}{a_2k_1}\leq \alpha,\ \ \
0<\frac{\gamma_{11}}{\eta_{11}}\leq \frac{\gamma_{12}}{\eta_{12}}
    =\frac{a_{1}k_{2}(M_{1}/k_{2}+b_{1})}{a_{2}k_{1}(M_{1}/k_{1}+b_{1})}\leq \alpha.
\end{eqnarray}
By the similar way used in the case {\bf (E2')},
we obtain that this theorem holds in the case {\bf (E3')}.
Therefore, the proof is now complete.
\hfill $\Box$
\subsection{Proof of Theorem \ref{thm-coexist}}
To complete the proof for Theorem \ref{thm-coexist},
we first consider the relation between the solutions of systems (\ref{2D-PP-model}) and (\ref{iterat-1}).
\begin{lemma}
\label{prpty-it-1}
Assume that either {\bf (A1)} or {\bf (A3)} holds.
Then for each $(U_{0},V_{0})\in\mathcal{A}$,
the solutions $(U,V)$ of system (\ref{2D-PP-model})
and $(U_{1j},V_{1j})$ of systems (\ref{iterat-1}) with $j=1,2$,
satisfy the following statements:
\begin{enumerate}
\item[{\bf (i)}]
$
U_{11}(t,x)\leq U(t,x)\leq U_{12}(t,x),\ \
V_{11}(t,x)\leq V(t,x)\leq V_{12}(t,x), \ \ t>0, \ x\in \Omega.
$
\item[{\bf (ii)}]
let $\widetilde{r}_{1j}$ be defined by
$\widetilde{r}_{1j}(x)=r_{1}(x)-\gamma_{1j}\theta_{d_{2},r_{2}}(x)/\eta_{1j}$ for $x\in \Omega$,
then
\begin{eqnarray}
\label{limt-U-V-1j}
\lim_{t\to+\infty}(U_{1j}(t,x),V_{1j}(t,x))=(\theta_{d_1,\widetilde{r}_{1j}}(x)/b_{1}, \theta_{d_2,r_2}(x)/\eta_{1j}), \ \ \ j=1,2,
\end{eqnarray}
and the following inequalities hold:
\begin{eqnarray}
&\theta_{d_1,\widetilde{r}_{11}}(x)/b_{1}\leq \theta_{d_1,\widetilde{r}_{12}}(x)/b_{1}, \ \ \ \
\theta_{d_2,r_2}(x)/\eta_{11}\leq \theta_{d_2,r_2}(x)/\eta_{12},\label{est-5}\\
& \theta_{d_1,\widetilde{r}_{11}}(x)/b_{1} \leq \liminf_{t\to+\infty} U(t,x)
\leq \limsup_{t\to+\infty}U(t,x)\leq \theta_{d_1,\widetilde{r}_{12}}(x)/b_{1},\nonumber\\
&\theta_{d_2,r_2}(x)/\eta_{11}
\leq \liminf_{t\to+\infty} V(t,x)
\leq \limsup_{t\to+\infty}V(t,x)
\leq \theta_{d_2,r_2}(x)/\eta_{12}.\nonumber
\end{eqnarray}
\end{enumerate}

\end{lemma}
{\bf Proof.}
Similarly to the case {\bf (E1')} in the proof of Theorem \ref{thm-global-1},
we obtain that {\bf (i)} holds.

To prove {\bf (ii)}, we only assume that the conditions stated in {\bf (A1)} hold,
the other case can be similarly discussed.
By the similar method used in the proof of Theorem \ref{thm-global-1}
we obtain that $\gamma_{1j}/\eta_{1j}\in (0,\alpha]$.
This together with Lemma \ref{lm-D-vo},
the fact that $\widetilde{\mathcal{D}}_{0}=\mathcal{D}_{0}$ for $a_{1}k_{2}/(a_{2}k_{1})=c_{1}/c_{2}$
and Lemma \ref{thm-stab-contrl} yields that
\begin{eqnarray*}
\widetilde{\mathcal{D}}_{+}=\varnothing, \ \ \ 
\widetilde{\mathcal{D}}_{0}=\left\{(d_1,d_2)\in\mathbb{R}_{+}^2 \ | \ r_1(x)\equiv \alpha\theta_{d_2,r_2}\right\}.
\end{eqnarray*}
Since $(d_1,d_2)\in\mathbb{R}_{+}^2\setminus\left\{(d_1,d_2)\in\mathbb{R}_{+}^2 \ | \ r_1(x)\equiv \alpha\theta_{d_2,r_2}\right\}$
in {\bf (A1)},
then $$(d_{1},d_{2})\in (\mathbb{R}_{+}^{2}\setminus (\widetilde{\mathcal{D}}_{+}\cup \widetilde{\mathcal{D}}_{0})).$$
Thus by Lemma \ref{prop-stab-coexist} {\bf (ii)} the limits in (\ref{limt-U-V-1j}) hold.
The left statements in {\bf (ii)} can be obtained by applying (\ref{est-4}) and {\bf (i)} in this lemma.
Therefore, the proof is now complete.
\hfill$\Box$

We next show the limit behavior of $U(t,x)$ as $t$ tends to infinity under the assumption that
either {\bf (A2)} or {\bf (A3)} holds.
More precisely, we have the following lemma.
\begin{lemma}
\label{lm-lim-U}
Assume that either {\bf (A2)} or {\bf (A3)} holds.
Then
there exists a unique positive  function $U^{*}(x)$ for $x\in \Omega$,
which is independent of $(U_{0},V_{0})$, such that
the solution $(U,V)$ of the modified Leslie-Gower model (\ref{2D-PP-model}) satisfies that
\begin{eqnarray*}
\lim_{t\to +\infty}U(t,x)=U^{*}(x), \ \ \ x\in \Omega.
\end{eqnarray*}
\end{lemma}
{\bf Proof.}
Since either {\bf (A2)} or {\bf (A3)} holds,
then by Lemma \ref{prpty-it-1} we have that
\begin{eqnarray}
\label{est-6}
\theta_{d_1,\widetilde{r}_{11}}(x)/b_{1} \leq \liminf_{t\to+\infty} U(t,x)
\leq \limsup_{t\to+\infty}U(t,x)\leq \theta_{d_1,\widetilde{r}_{12}}(x)/b_{1}.
\end{eqnarray}
Since $k_1=k_2$, then
\begin{eqnarray*}
\frac{\gamma_{11}}{\eta_{11}}=\frac{\gamma_{12}}{\eta_{12}}=\frac{a_{1}}{a_{2}},
\end{eqnarray*}
which yields that
\begin{eqnarray*}
\widetilde{r}_{11}(x)=r_{1}(x)-\frac{\gamma_{11}}{\eta_{11}}\theta_{d_{2},r_{2}}(x)
=r_{1}(x)-\frac{\gamma_{12}}{\eta_{12}}\theta_{d_{2},r_{2}}(x)=\widetilde{r}_{12}(x).
\end{eqnarray*}
Hence, by the above equalities and (\ref{est-6}) the proof is finished.
\hfill$\Box$

In the end of this section we give the proof for Theorem \ref{thm-coexist}.\\
{\bf Proof of Theorem \ref{thm-coexist}}.
If either {\bf (A1)} or {\bf (A3)} holds,
then by Lemma \ref{prpty-it-1} we have that
\begin{eqnarray*}
\liminf_{t\to+\infty}\min_{x\in \overline{\Omega}} U(t,x)>\delta>0,\ \ \
\liminf_{t\to+\infty}\min_{x\in \overline{\Omega}} V(t,x)>\delta>0,
\end{eqnarray*}
where $\delta$ satisfies that
$$2\delta=\min\left\{\min_{x\in \overline{\Omega}}\theta_{d_1,\widetilde{r}_{11}}(x)/b_{1},\,
\min_{x\in \overline{\Omega}}\theta_{d_2,r_2}(x)/\eta_{11}\right\}.$$
Then {\bf (i)} is proved.

Assume that system (\ref{2D-PP-model}) has a positive coexistence steady state $(u^{*},v^{*})$,
then $u^{*}=U^{*}$, where $U^{*}$ is a positive function obtained as in Lemma \ref{lm-lim-U}.
To obtain $v^{*}$ we consider the elliptic equation in the form
\begin{equation*}
\left\{
\begin{aligned}[cll]
&d_{2}\Delta V+V\left(r_{2}(x)-\frac{a_{2}V}{U^{*}(x)+k_{2}}\right)=0\  \ \ \  \ && \mbox{in}\ \Omega\times\mathbb{R}_{+},\\
&\frac{\partial V}{\partial n}=0 \  \ \ \  \ & &  \mbox{on}\ \partial\Omega\times\mathbb{R}_{+},\\
\end{aligned}
\right.
\end{equation*}
which is equivalent to
\begin{equation}
\label{ellp}
\left\{
\begin{aligned}[cll]
&d_{2}\Delta V+\frac{V}{U^{*}(x)+k_{2}}\left(r_{2}(x)(U^{*}(x)+k_{2})-a_{2}V\right)=0\  \ \ \  \ && \mbox{in}\ \Omega\times\mathbb{R}_{+},\\
&\frac{\partial V}{\partial n}=0 \  \ \ \  \ & &  \mbox{on}\ \partial\Omega\times\mathbb{R}_{+},\\
\end{aligned}
\right.
\end{equation}
Let the constants $W_{*}$ and $W^{*}$ be given by
\begin{eqnarray*}
W_{*}=\min_{x\in\overline{\Omega}}\frac{r_{2}(x)(U^{*}(x)+k_{2})}{a_{2}},\ \ \ \
W^{*}=\max_{x\in\overline{\Omega}}\frac{r_{2}(x)(U^{*}(x)+k_{2})}{a_{2}}.
\end{eqnarray*}
Since $U^{*}$ is a positive function and $r_{2}$ is a non-negative function with $\overline{r}_{2}>0$,
then $W_{*}\geq 0$ and $W^{*}>0$,
which implies that  $W_{*}$ is a lower solution of (\ref{ellp})
and $W^{*}>0$ is a upper-solution.
Consequently, by the method of upper and lower solutions,
we have that (\ref{ellp}) has a unique positive solution $v^{*}$, which is globally asymptotically stable.
This together with Lemma \ref{lm-lim-U} yields that {\bf (ii)} holds.
Therefore, the proof is finished.
\hfill$\Box$

\section{Concluding remarks}
\setcounter{equation}{0}
\setcounter{lemma}{0}
\setcounter{theorem}{0}
\setcounter{remark}{0}

In this paper we have investigated the global dynamics of the  competition-diffusion system (\ref{control equation}),
which indicates that the evolution of the density of the predator is independent of the density of the prey.
By the  principal spectral theory and the dynamics of the classical single-species logistic model,
we prove that its global dynamics can be determined by the local stability of its semi-trivial steady state $(0,\theta_{d_{2},r_{2}}/c_{2})$.
We also apply the obtained results to 
obtain the global stability of the steady states and the persistence of the modified Leslie-Gower model (\ref{2D-PP-model})
under some suitable conditions.
It is interesting to give a more detailed study of the positive coexistence steady states for system (\ref{2D-PP-model})
satisfying the following condition:
\begin{eqnarray*}
\frac{a_1k_2}{a_2k_1}\in (0,\alpha],\ \ k_2>k_1,\ \
(d_1,d_2)\in\mathbb{R}_{+}^2\setminus\left\{(d_1,d_2)\in\mathbb{R}_{+}^2 \ | \ r_1(x)\equiv \alpha\theta_{d_2,r_2}\right\}.
\end{eqnarray*}

\section*{Acknowledgments}

This work was partly supported by the Hubei provincial postdoctoral science and technology activity project (No. Z2).

\end{document}